\documentstyle[leqno,12pt]{amsart}
\def\ds{\displaystyle}

\newsymbol\Subset 1362
\def\ds{\displaystyle}

\newcommand{\R}{{\Bbb R}}

\newcommand{\M}{{\Bbb M}}
\newcommand{\A}{{\Bbb A}}

\newcommand{\loc}{{\rm loc}}

\newcommand{\const}{{\rm const}}

\def\Xint#1{\mathchoice 
{\XXint\displaystyle\textstyle{#1}}% 
{\XXint\textstyle\scriptstyle{#1}}% 
{\XXint\scriptstyle\scriptscriptstyle{#1}}% 
{\XXint\scriptscriptstyle\scriptscriptstyle{#1}}% 
\!\int} 
\def\XXint#1#2#3{{\setbox0=\hbox{$#1{#2#3}{\int}$ } 
\vcenter{\hbox{$#2#3$ }}\kern-.6\wd0}}

\def\ol{\overline}%%%%%%%

\def\B{{\mathcal B}}
\def\C{{\mathcal C}}

\def\F{{\Bbb  F}}

\def\O{{\mathcal O}}

\def\bF{{\mathbf F}}
\def\bA{{\mathbf A}}

\def\div{{\rm div}}
\def\const{{\rm const}}
\def\loc{{\rm loc}}
\def\ba{{\mathbf a}}
\def\bu{{\mathbf u}}
\def\bz{{\mathbf z}}
\def\bg{{\mathbf g}}
\def\ff{{\mathbf f}}
\def\bb{{\mathbf b}}

\newtheorem{thm}{Theorem}[subsection]  %% Definition of Theorem
\newtheorem{defin}[thm]{Definition}    %% Definition of Definition
\makeatletter
\def\theequation{\@arabic\c@equation}
\def\thethm{\@arabic\c@thm}
\def\thelem{\@arabic\c@thm}
\def\thecrlr{\@arabic\c@thm}
\def\theprp{\@arabic\c@thm}
\def\therem{\@arabic\c@thm}
\makeatother

\begin{document}

\baselineskip=16pt
%\today

\title[Quasilinear  systems in divergence form] {$L^p$-integrability of the gradient of  solutions to  quasilinear  systems with discontinuous coefficients}

\author[L.G. Softova]{Lubomira G. Softova}

\address{Department of Civil Engineering\\
  Second University  of Naples\\
 Via Roma 29\\
 81031 Aversa\\
 Italy}

\email{luba.softova@@unina2.it}



\maketitle

\begin{abstract}
The Dirichlet problem for a class of quasilinear elliptic systems of equations with 
small-BMO coefficients in Reifenberg-flat domain $\Omega$ is considered.
 The lower order  terms  supposed to satisfy controlled growth conditions in $\bu$ and $D\bu.$  It is obtained $L^p$-integrability  with $p>2$ of  $D\bu$  where $p$ depends explicitly on the  data. An analogous result is obtained also for the Cauchy-Dirichlet problem for  quasilinear parabolic systems.
\end{abstract}

\section{Introduction}
In the present work we study the integrability properties of the weak solutions of the following Dirichlet problem
\begin{equation}\label{DP}
\begin{cases}
 D_\alpha\big(A^{\alpha\beta}_{ij}(x) D_\beta u^{j} (x)
+a^{\alpha}_i(x,\bu) \big)= b_i(x,\bu,D\bu) & \text{ a.a. } x\in\Omega \\
\bu(x)=0\quad  & \text{ on } \partial\Omega
\end{cases}
\end{equation}
where $\Omega\subset\R^n,$ $n\geq 2$ is a bounded {\it Reifenberg flat} domain (see Definition~\ref{reifenberg}).  The principal coefficients are  discontinuous with 
"small" discontinuity expressed in terms of their   {\it bounded mean oscillation} (BMO) in $\Omega$
 (cf. \cite {JN}).
The matrix 
$\bA(x)=\{A_{ij}^{\alpha\beta}(x)\}_{i,j\leq N}^{\alpha,\beta\leq n}$ verifies 
\begin{equation}\label{eq4}
A^{\alpha\beta}_{ij}(x)\xi^i_\alpha\xi^j_\beta\geq \lambda |\xi|^2 \quad\forall \ \xi\in \M^{N\times n},\quad
\|\bA\|_{\infty,\Omega}\leq M
\end{equation} 
with some positive constants $\lambda$ and $M.$  
The non linear terms 
$$
\ba(x,\bu)=\{a^\alpha_i(x,\bu)\}^{\alpha\leq n}_{i\leq N} \quad \text{ and  } \quad 
\bb(x,\bu,\bz)=\{b_i(x,\bu,\bz)\}_{i\leq N}
$$
 supposed to be  Carath\'eodory functions for
 $x\in\Omega,$  $\bu\in \R^N,$ $\bz \in \M^{N\times n}$ and satisfy {\it controlled growth conditions}. Namely,  for  $  |\bu|, |\bz|\to\infty$ we have 
 \begin{equation}\label{eq0}
a^\alpha_i(x,\bu)=\O(\varphi_1(x)+|\bu|^{\frac{n}{n+2}}) \quad \text{ and }
\quad b_i(x,\bu,\bz)=\O(\varphi_2(x)+|\bu|^{\frac{n+2}{n-2}}) 
\end{equation}
 with $\varphi_1\in L^p(\Omega),$ $p>2$ and $\varphi_2\in L^q(\Omega),$ $q>\frac{2n}{n+2}.$

Our aim is to show that the problem  \eqref{DP} satisfies 
  the   Calder\'on--Zygmund property when $\Omega$ is $(\delta,R)$-Reifenberg flat  and the coefficients are $(\delta,R)$-vanishing in $\Omega.$ Precisely, each  bounded
 weak solution $\bu\in W_0^{1,2}\cap L^\infty(\Omega;\R^N)$ of \eqref{DP} gains better regularity from the data $\varphi_1$ and $\varphi_2$ and belongs to $W_0^{1,\min\{p,q^\ast\}}\cap L^\infty(\Omega;\R^N)$ where $q^\ast$ is the Sobolev conjugate of $q$  (see  \eqref{conjugate}).

Similar result is obtained also for the Cauchy-Dirichlet problem for  the parabolic  quasilinear system 
$$
\begin{cases}
u_t^i-D_\alpha(A^{\alpha\beta}_{ij}(x,t)D_\beta u^j+a^\alpha_i(x,t,\bu))=b_i(x,t,\bu,D\bu)   ) & \text{ a.a. }  (x,t)\in Q\\
\bu(x,t)=0 & (x,t)\in \partial Q
\end{cases}
$$
in a cylinder $Q=\Omega\times(0,T)$ where $\Omega $ is $(\delta,R)$-Reifenberg flat and
 $\partial Q=\Omega\cup \{\partial\Omega\times(0,T)\}$ is the parabolic boundary.

The  problem of integrability and regularity  of the solutions of linear and  quasilinear   elliptic/parabolic equations and systems is widely  studied.  Let us start with the classical results concerning equations/systems with smooth coefficients  presented in the monographs \cite{LU, LSU}.   In the scalar case, $N=1,$  the notorious  results of De Giorgi \cite{DeG} and Nash \cite{N} assert H\"older regularity of the solutions of linear   divergence form equations with only $L^\infty$ principal coefficients. One remarkable result that permits to obtain higher integrability of the weak solutions is due to Gehring \cite{Gh}. He studied  
 integrability properties of functions satisfying the reverse H\"older inequality. It was noticed that some power of the gradient of the weak solutions satisfies  local reverse H\"older inequality.  
Modifying the Gehring lemma,  Giaquinta and Modica \cite{GM} firstly obtain higher integrability  of solutions of divergence form quasilinear elliptic  equations. For the sake of completeness we give this result as it is presented in the monograph by Giaquinta \cite[Theorem~V.2.3]{G}.
\begin{thm}\label{thmG}
Suppose that $g\in L^q(\Omega),$ $F\in L^{q+\delta}(\Omega),$ $g,F\geq 0,$ $q>1,$ $\delta>0$ and 
$$
\Xint-_{\B_R(x)}g^qdx\leq B\left( \Xint-_{\B_{2R}(x)}gdx\right)^q +\Xint-_{\B_{2R}(x)}F^qdx+
\theta \Xint-_{\B_{2R}(x)} g^qdx
$$
for a.a. $x\in \Omega,$ $R<\frac12\min\{ d(x,\partial\Omega), R_0 \}$ where  $R_0>0,$ $B>1,$ $\theta\in[0,1).$ Then $g\in L^{p,\loc}(\Omega)$ and 
$$
\left( \Xint-_{\B_R(x)}g^pdx \right)^{1/p}\leq C\left\{ \left( \Xint-_{\B_{2R}(x)}g^qdx \right)^{1/q}
+ \left( \Xint-_{\B_{2R}(x)}F^pdx \right)^{1/p}   \right\}
$$
for any  ball $\B_{2R}\subset \Omega,$ $2R<R_0,$ $p\in[q,p_0)$ where $C>0,$ $p_0>q$ depend only on $B,\theta,q,n.$
\end{thm}
There are various generalizations of the above  theorem permitting to study elliptic and parabolic problems with Dirichlet and Neumann boundary conditions (see \cite{A,A2,A3,AL,DK,FZ,MM}).
Other results concerning higher integrability  of divergence form  quasilinear equations and variational equations  could be found in  \cite{C,G, GS, JS, MPS}.   
 The $L^p$-estimates of derivatives obtained such way laid the foundation to the so-called "direct method" of proving partial regularity of solutions.  Recently, the method of A-harmonic approximation permits to study the regularity of the solutions without the use of the Gehring lemma. For more details we refer the reader to  \cite{AM,DG,DKM,DMS}, see also the references therein.  

The regularity theory for linear  operators  with smooth data   was extended on operators with discontinuous coefficients defined in rough domains.  In  \cite{B1, BR, BW}      the authors consider  divergence form elliptic and parabolic  equations and systems  with $BMO$ coefficients in Reifenberg flat domain with Dirichlet boundary conditions  extending  such way   the known results on  operators with $VMO$ coefficients too (see also \cite{MPS,PSpsyst,PSesyst,K1,K2} and the references therein).   In \cite{FZ}  a reverse H\"older inequality is established for quasilinear elliptic systems with principal coefficient being  $VMO$ in $x$ and under controlled growth conditions over the lower order terms.  It permits the authors to obtain interior H\"older continuity of solutions to scalar equations as well as partial H\"older regularity of solutions to systems. In \cite{P0,P}  global H\"older regularity of solutions to elliptic quasilinear equations with $VMO$ in $x$ principal coefficients is proved under {\it strictly controlled growth} conditions. Later this result is extended for quasilinear elliptic and parabolic equations in Reifenberg flat domains supposing {\it controlled growth} conditions and Dirichlet boundary data (see \cite{DK,P1,PS1,PS2}).

In the present work  we extend the results from \cite{PS2} to elliptic and parabolic systems with discontinuous data. Making use of the linear  $L^p$-theory for  systems, developed  in \cite{BR,BW}  and    the bootstrap method we prove $D\bu\in L^r$
with  $r$ depending explicitly  on the data $\varphi_1$ and $\varphi_2$ in \eqref{eq0}.

\section{Elliptic systems, definitions and main result}

 In the following we use the  standard  notations:\\
$\bullet$ $x=(x_1,\ldots,x_n)\in \R^n,$ $\rho>0$  and  $\B_\rho(x) = \{y\in \R^n: |x-y|<\rho\}.$\\
$\bullet$ let  $\Omega\subset\R^n$ be a bounded domain,  $x\in \Omega$ and denote  
 $\Omega_\rho(x)=\Omega\cap \B_\rho(x).$\\
$\bullet$  $\M^{N\times n}$  is   the set of  $N\times n$-matrices.\\
$\bullet$ For a vector function  $\bu= (u^1,\ldots,u^N):\Omega\to \R^N $ we  write
$$
|\bu|^2=\sum_{j\leq N} |u^j|^2, \quad   D_\alpha u^j=\frac{\partial }{\partial x_\alpha}u^j, 
$$
$$
  D\bu=\{D_\alpha u^j\}^{\alpha\leq n}_{j\leq N}\in \M^{N\times n}, \quad  
 |D\bu|^2=\sum_{\underset{j\leq N}{\alpha \leq n}}|D_\alpha u^j|^2 .
$$
$\bullet$ Let $f:\Omega\to \R$ and  $|\Omega|$  be the Lebesgue measure of $\Omega,$ then  
$$
\Xint-_\Omega f(y) dy=\frac1{|\Omega|}\int_\Omega f(y) dy,\quad 
 \|f\|^p_{p,\Omega}=\|f\|^p_{L^p(\Omega)} =\int_\Omega |f(y)|^pdy.
$$
$\bullet$ For $\bu\in L^p(\Omega;\R^N)$ write $\|\bu\|_{p,\Omega}$ instead of
 $\|\bu\|_{ L^p(\Omega;\R^N)}.$\\
$\bullet$ 
For each   $s\in(1,\infty)$  recall   that  $s^\ast$  means the Sobolev conjugate of $s$ 
\begin{equation}\label{conjugate}
s^\ast=\begin{cases}
\frac{ns}{n-s} & \text{  if } s<n\\
\text{ arbitrary large number } > 1 & \text{  if  } s\geq n.
\end{cases}
\end{equation}
For the function spaces we follow the notions of the monographs \cite{LU, MPS}.
 Through all the paper the standard summation convention on repeated upper and lower indexes is adopted. 
The letter $C$ is used for various constants and may change from one occurrence to another.

In \cite{R}  Reifenberg  introduced a class of domains with rough  boundary  that  can be approximated by hyperplanes at every point and at every scale. Namely 
\begin{defin}\label{reifenberg}\rm
The domain $\Omega$ is $(\delta,R)$-Reifenberg flat if there exist positive constants $R,$ $\delta<1 $ such that for each $x\in \partial \Omega$
and each $\rho\in(0,R)$ there is a local coordinate system $\{y_1,\ldots,y_n\}$ with the property
\begin{equation}\label{eq2}
\B_\rho(x)\cap \{y_n>\delta\rho \}\subset\Omega_\rho(x)\subset\B_\rho(x)\cap\{y_n>-\delta\rho \}.
\end{equation}
\end{defin}
Reifenberg  arrived at that concept of flatness in his studies on Plateau's problem  in higher dimensions and he proved that such a domain is locally a topological disc when $\delta$  is small enough.   It is easy to see that 
a  $C^1$-domain is  a  Reifenberg flat  with $\delta\to 0 $ as $R\to 0.$  A domain with Lipschitz boundary with a Lipschitz  constant less than $\delta$ also verifies the condition \eqref{eq2} if  $\delta$  is small enough (say $\delta<1/8$).     But the class of Reifenberg's  domains is much  more wider  and contain domains with fractal boundaries. For instance, consider
a self-similar snowflake  $S_\beta.$  It is  a flat version of the Koch snowflake $S_{\pi/3}$ where the angle of the spike with respect to the horizontal is $\beta.$    A domain $\Omega\subset\R^2$ with $S_\beta\subset\partial\Omega$ is a Reifenberg flat  if $0<\sin\beta<\delta<1/8.$
This kind of flatness exhibits  minimal geometrical conditions necessary for some natural properties in analysis and potential theory to hold.  For more detailed overview  of the properties of these  domains
we refer the reader to the papers \cite{MT,T}.

From \eqref{eq2} it follows that  $\partial\Omega$ satisfies the  $(A)$-property  (cf. \cite{C, G, LU}).  
Precisely, the  measure $|\Omega_\rho(x)|$ is $\delta$-comparable to  $|\B_\rho(x)|,$ that is there exists a positive  constant  $A(\delta)<1/2$ such that 
\begin{equation}\tag{A}\label{A}
 A(\delta)|\B_\rho(x)| \leq |\Omega_\rho(x)|\leq (1- A(\delta))|\B_\rho(x)|
\end{equation}
for any fixed $x\in\partial\Omega,$ $\rho\in(0,R)$ and $\delta\in(0,1).$  This condition excludes that $\Omega $  may have sharp outward and inward cusps.
Moreover, for small $\delta$ they can be approximated  in a uniform way by Lipschitz domains with a Lipschitz constant less then $\delta$  (see \cite[Lemma~5.1]{BW}). As consequence, they  are $W^{1,p}$-extension domains, $1\leq p\leq \infty,$ hence 
the usual extension theorems, the  Sobolev and Sobolev--Poincar\'e inequalities are valid in $\Omega.$

To describe the discontinuity of the principal  coefficients we need of the following
\begin{defin}\label{def2}
We say that a function $a(x)$ is a $(\delta,R)$-vanishing if there exist positive constants $R$ and $ \delta<1$ such that 
\begin{align}\label{eq5}
\sup_{0<\rho\leq R}\sup_{x\in\Omega}\Xint-_{\Omega_\rho(x)}|a(y) -\ol a_{ \Omega_\rho(x)}   |^2 dy\leq \delta^2,\quad
\ol a_{ \Omega_\rho(x)}= \Xint-_{\Omega_\rho(x)} a(y) dy.
\end{align}
\end{defin}
We suppose that all $A^{\alpha\beta}_{ij}(x)$ are $(\delta,R)$-vanishing. It implies   that $\bA\in BMO(\Omega)$ with a small BMO norm $\|\bA\|_\ast<\delta.$ 

The nonlinear terms 
    $\ba(x,\bu)$ and $\bb(x,\bu,\bz)$ are Carath\'eodory functions for $x\in\Omega,$  $\bu\in \R^N,$ $\bz \in \M^{N\times n}$ and satisfy the {\it controlled growth conditions} 
\begin{align}\label{contr1}
&  |\ba(x,\bu) |\leq \Lambda(\varphi_1(x)+|\bu|^{\frac{n}{n-2}}), \quad    \varphi_1\in L^p(\Omega), \ p>2 \\
\label{contr2}
&|\bb(x,\bu,\bz)|\leq \Lambda\big(\varphi_2(x) + |\bu|^{\frac{n+2}{n-2}}+|\bz|^{\frac{n+2}{n}}  \big), \quad \varphi_2\in L^q(\Omega), \ q>\frac{2n}{n+2}
\end{align}
with some positive constant $\Lambda.$
In the particular case  $n=2$ the powers of $|\bu|$ could be arbitrary positive numbers while the growth of $|\bz|$ is quadratic (cf. \cite{G, LU}).

Under a {\it weak solution} to the problem \eqref{DP} we mean a function  
 $\bu\in W_0^{1,p}(\Omega;\R^N),$ $1<p<\infty$  satisfying
\begin{align*}
\int_\Omega& A^{\alpha\beta}_{ij}(x)D_\beta u^j(x)D_\alpha\chi^i(x) dx +\sum_{\underset{i\leq N}{\alpha\leq n}}\int_\Omega 
a^\alpha_i(x,\bu(x))D_\alpha\chi^i(x)dx\\
& +\int_\Omega b_i(x,\bu(x),D \bu(x))\chi^i(x) dx=0,\quad j=1,\ldots,N
\end{align*}
for all $\chi\in W_0^{1,p'}(\Omega;\R^N),$ $p'=p/(p-1).$ 
The conditions  \eqref{contr1}  and \eqref{contr2}
are the natural ones that ensure convergence of the integrals above. Moreover, they are optimal since a growth of  the gradient  greater than $\frac{n+2}n$ leads to unbounded solutions as it is seen from the following example (cf. \cite{LSU, P}). The function $u(x)\in W^{1,2}(\B_1(0)),$ $u(x)=|x|^{\frac{r-2}{r-1}}$  is a solution  of the equation $\Delta u=C|Du|^r$  in $\B_1(0).$ Note that $u(x)\not\in L^\infty(\B_1(0))$ for $\frac{n+2}{n}<r<2.$

In generally we cannot expect boundedness of each solution of \eqref{DP} unless we add some structural conditions. Consider, for instance, the system
$$
D_\alpha(A^\alpha_i(x,\bu,D\bu))=b_i(x,\bu,D\bu)\qquad x\in \Omega
$$
where 
$$
A_i^\alpha(x,\bu,D\bu)=\sum_{\beta\leq n}\sum_{j\leq N}(A^{\alpha\beta}_{ij}(x) D_\beta u^j +a^\alpha_i(x,\bu) )
$$
are measurable in $x\in\Omega.$  Assume a pointwise coercive and sign conditions, both of them for large values of the corresponding component of $\bu,$ precisely:  for every $i\in\{1,\ldots,N\}$ there exist constants  $\theta^i,M^i,\nu \in(0,+\infty)$ such that  for $u^i\geq \theta^i$ we have
\begin{equation}\label{c1}
\begin{cases}
\ds \nu|\xi^i|^2-M^i\leq \sum_{\alpha\leq n}A_i^\alpha (x,\bu,\xi)\xi^\alpha_i\\
  \ds 0\leq b_i(x,\bu,\xi)\quad  \text{ for a.a. } x\in\Omega, \  \forall \xi\in \M^{N\times n}.
\end{cases}
\end{equation}
  Suppose  \eqref{contr1}, \eqref{contr2} and  \eqref{c1}  and let $\bu\in W^{1,2}\cap L^{\frac{2n}{n-2}}(\Omega;\R^N)$ be a weak solution of \eqref{DP} then for each $i\in \{1,\ldots,N\}$
$$
\sup_\Omega u^i\leq \theta^i + K^i
$$
where $K^i$ depend on $M^i, n, |\Omega|$ and $ \nu$ (see \cite{LP}).

 \begin{thm}\label{th1}
Let $\bu\in W_0^{1,2}\cap L^\infty(\Omega;\R^N)$ be a weak solution of the problem \eqref{DP} under the conditions 
\eqref{eq4}, \eqref{contr1} and  \eqref{contr2}.
Then there exists a small number $\delta_0>0$ such that if 
$\Omega$ is $(\delta,R)$-Reifenberg flat domain and $A^{\alpha\beta}_{ij}(x)$ are $(\delta,R)$-vanishing with $\delta<\delta_0<1$ then
\begin{equation}\label{qast}
\bu\in W_0^{1,r}\cap  L^\infty(\Omega;\R^N)\  \text{ with }\  r=\min\{p,q^\ast\}. 
\end{equation}
\end{thm}
\begin{pf}
In \cite[Chapter~5]{G} Giaquinta considers quasilinear strongly  elliptic systems with $L^\infty$ principal coefficients,  under the conditions \eqref{contr1} and \eqref{contr2}.  Making use of the reverse  H\"older's inequality and the version of the  Gehring lemma  it is shown that  there exists an exponent $r_0>2 $ such that $\bu\in W^{1,r_0}_{\loc}(\Omega,\R^N)$
 (cf. \cite[Theorem~V.2.3]{G},\cite[Chapter~III]{C} or \cite[Lemma~3.2.23]{MPS} ).  Since, roughly speaking, Caccioppoli-type inequalities hold up to the boundary, the method for obtaining higher integrability can be carried over up to the boundary. In \cite[Chapter~5]{G} it is done for the Dirichlet problem in Lipschitz domain. Since the Reifenberg flat domain 
can be uniformly approximated by Lipschitz domains the same result still holds true. Precisely, there is $r_0>2$ such that
    \begin{equation}\label{eq1}
\|D\bu\|_{r,\Omega}\leq N \qquad \forall\  r\in[2,r_0)
\end{equation}
where $N$ and $r_0$ depend on  $n, \Lambda, \lambda, \|\varphi_1\|_{p,\Omega},
 \|\varphi_2\|_{q,\Omega}, 
 |\Omega|, \|D\bu\|_{2,\Omega}.$

Let  $n>2$  and   $\bu\in W_0^{1,2}\cap L^\infty(\Omega;\R^N)$ be a solution of  \eqref{DP}. Fixing that solution in the nonlinear terms we get the linearized problem
\begin{equation}\label{eq3}
\begin{cases}
\ds D_{\alpha}(A^{\alpha\beta}_{ij}(x)D_\beta u^j )=f_i(x)-\div (\A_i(x)) & \text{ a.a. }  x\in\Omega\\
\bu(x)=0  & \text{ on }  \partial\Omega
\end{cases}
\end{equation}
where
\begin{align*}
f_i(x)=& b_i(x,\bu,D\bu),\qquad \ff(x)=\bb(x,\bu,D\bu),\\
\A_i(x)=& (a^1_i(x,\bu),\ldots, a_i^n(x,\bu)), \qquad  \A(x)=(\A_1(x),\ldots,\A_N(x))
\end{align*}
and by  \eqref{contr1}, \eqref{contr2} and \eqref{eq1} we get
 \begin{equation}\label{estimates}
\begin{cases}
\ds \|\A\|_{p,\Omega}\leq C \left(\|\varphi_1\|_{p,\Omega} +\|\bu\|^{\frac{n}{n+2}}_{\infty,\Omega} \right)\\[6pt]
\ds \|\ff\|_{q_1,\Omega}\leq C\left( \|\varphi_2\|_{q_1,\Omega} +\|\bu\|^{\frac{n+2}{n-2}}_{\infty,\Omega} +
 \|D\bu\|_{\frac{q_1(n+2)}{n},\Omega}^{\frac{n+2}{n}}   \right) 
\end{cases}
\end{equation}
with $p>2$ and $  q_1=\min\big\{q, \frac{r_0n}{n+2}\big\}.$
 Further, for all $f_i\in L^{q_1}(\Omega),$ $i=1,\ldots,N$
there exists a vector field 
$\F_i(x)\in L^{q_1^\ast}(\Omega,\R^n)$ such that $f_i(x)=\div\, \F_i(x).$ Denote $\F(x)=(\F_1(x),\ldots, \F_N(x)),$ then by   \cite[Lemma~3.1]{P}) we have 
\begin{equation}\label{q1}
\|\F\|_{q_1^\ast,\Omega}\leq C\|\ff\|_{q_1,\Omega},\quad  q_1=\min\big\{q,\frac{r_0n}{n+2}\big\}\,.
\end{equation}
  Thus the problem \eqref{eq3} becomes
\begin{equation}\label{linearized}
\begin{cases}
\ds  D_\alpha(A^{\alpha\beta}_{ij}(x)D_{\beta}u^j(x) ) =\div(\F_i(x)-\A_i(x)) & \text{ a.a. } x\in \Omega\\
\bu(x)=0, & x\in \partial\Omega.
\end{cases}
\end{equation}
For linear  systems as above  we dispose with the regularity result of Byun and Wang \cite[Theorem~1.7]{BW} that asserts there
 exists a small positive constant 
$\delta=\delta(\lambda, p,n,N)$ such that for each $(\delta,R)$-vanishing
 $A^{\alpha\beta}_{ij},$ for each $(\delta,R)$-Reifenberg flat $\Omega,$ and for each matrix function 
$\F-\A\in L^{r_1}(\Omega;\M^{N\times n}),$ 
with $r_1=\min\{p,q_1^\ast\},$  
 the solution 
$\bu\in W_0^{1,2}\cap L^\infty(\Omega;\R^N)$ of \eqref{linearized} belongs to  
$ W_0^{1,r_1}\cap L^\infty(\Omega;\R^N)$ and the following estimate holds
\begin{equation}\label{r1}
\|D\bu\|_{r_1,\Omega}\leq C \|\F- \A\|_{r_1,\Omega},\qquad r_1=\min\{p,q_1^\ast\}
\end{equation}
with $C=C(\lambda, p,n,N,|\Omega|).$

Our goal is to show the inclusion $D\bu\in L^r(\Omega;\M^{N\times n})$ with  $r= \min\{p,q^\ast\}.$ For this  we  study  the following cases:
\begin{itemize}
\item[$1)$] If   $q\leq \frac{r_0n}{n+2}$ then  $q_1=q$   in \eqref{q1}  and $r_1\equiv r=\min\{p,q^\ast\}.$
\item[$2)$]  If   $q>\frac{r_0n}{n+2},$ then $q_1= \frac{r_0n}{n+2}$  and
$$
q_1^\ast=\begin{cases}
\ds\frac{r_0n}{n+2-r_0} &\ds \text{ if  }  \    \frac{r_0n}{n+2}<n\\[10pt]
\ds\text{ arbitrary large number }> 1  &\ds  \text{ if } \   \frac{r_0n}{n+2}\geq n.
\end{cases}
$$
Consider again   two sub-cases:
\begin{itemize}
\item[$2_a)$] If   $n>\frac{r_0n}{n+2}$         then  $r_1=\min\{p,\frac{r_0n}{n+2-r_0}\}.$ 
If $r_1=p$ then the theorem holds true otherwise  $D\bu\in L^{r_1}(\Omega:\M^{n\times N})$ with  $r_1=\frac{r_0n}{n+2-r_0}.$
\item[$2_b)$] If   $n\leq \frac{r_0n}{n+2}$ then   $q_1^\ast$ is arbitrary large number, that implies  $r_1=p$  and the theorem holds true once again.
\end{itemize}
\end{itemize}
It is easy to see that $r_1\equiv r$ unless 
$$
\frac{r_0n}{n+2}<q\qquad \text{and}\qquad \frac{r_0n}{n+2}<n
$$
when  $\bu\in W^{1,r_1}_0\cap L^\infty(\Omega;\R^N)$ with $r_1=\frac{r_0n}{n+2-r_0}.$ It holds for any solution of the linearized problem \eqref{linearized} including the one fixed in the coefficients in \eqref{DP}. 
  
 Consider once  again \eqref{estimates} with  
$D\bu \in L^{r_1}(\Omega;\M^{N\times n}).$
Hence    $f_i\in L^{q_2}(\Omega)$  with 
$q_2=\min\{q, \frac{r_0n^2}{(n+2)^2-r_0(n+2)}\}$ and 
the associated vector-field  $\F_i$  belongs to $ L^{q_2^\ast}(\Omega;\R^n).$ Than 
 $\F_i-\A_i\in L^{r_2}(\Omega;\R^n)$ with $r_2=\min\{p,q_2^\ast\}.$
Applying   \cite[Theorem~1.7]{BW}
 to   system \eqref{linearized}
and repeating the same procedure as above we get  that the theorem holds with $r_2\equiv r$ if 
\begin{equation}\label{si}
i)\  \frac{r_0n^2}{(n+2)^2-r_0(n+2)}\geq n\qquad \text{ or }\qquad 
ii)\  
\frac{r_0n^2}{(n+2)^2 -r_0(n+2)}\geq q.
\end{equation}
Otherwise  $r_2=\frac{r_0n^2}{(n+2)^2-r_0(n+2)-r_0n}$ if  
\begin{equation}\label{no}
\frac{r_0n^2}{(n+2)^2-r_0(n+2)}<q\qquad \text{ and }\qquad \frac{r_0n^2}{(n+2)^2-r_0(n+2)}<n
\end{equation}
  Repeating the same procedure $k$-times we get that the assertion holds if 
\begin{equation}\label{eq5a}
\frac{r_0n^k}{(n+2)^k-r_0\sum_{s=0}^{k-2}n^s(n+2)^{k-1-s}}\geq \min\{n,q\}\,.
\end{equation}
 Direct calculations give that \eqref{eq5a} is equivalent to
$$
k>\min \left\{ \big[\log\frac{r_0}{r_0-2}\big/ \log\frac{n+2}{n} \big], 
\big[\log\frac{r_0(2q+qn+2)}{q(n+2)(r_0-2)}\big/\log\frac{n+2}{n}  \big]+1  \right\}
$$
 where $[x]$ means the integer part of $x.$

 The case $n=2$ is simpler   and is left to the reader.   
\end{pf}

\section{Quasilinear parabolic systems}

Let $Q=\Omega\times(0,T) $ be a cylinder in $\R^{n+1}$ with  $\Omega $ being    $(\delta,R)$-Reifenberg flat.
Denote by $\C_\rho$  the parabolic cylinder
$$
\C_\rho(x,t)=\B_\rho(x)\times(t-\rho^2,t),\quad Q_\rho(x,t)=Q\cap \C_\rho(x,t) \quad\text{ for } (x,t)\in Q,
$$
$$
\ol{a}_{Q_\rho(x,t)}=\Xint-_{Q_\rho(x,t)}a(y,\tau)dyd\tau=\frac1{|Q_\rho(x,t)|}\int_{Q_\rho(x,t)}
 a(y,\tau) dyd\tau.
$$
 Let $1<r<\infty$ and $\bu:Q\to \R^N.$\\
 1. The  space $W^{1,0}_r(Q;\R^N)$  consists of all  functions $\bu\in L^r(Q;\R^N)$ having a  finite norm
$$
\|\bu\|^r_{W^{1,0}_r(Q;\R^N)} = \|\bu\|^r_{r,Q} +\|D\bu\|^r_{r,Q}\,.
$$
 2. The space $W^{1,r}_\ast(Q;\R^N)= L^r(0,T; W^{1,r}(\Omega;\R^N))\cap W^{1,r}(0,T; W^{-1,r'}(\Omega;\R^N)) ,  $ $r'=r/(r-1)$  consists of the functions $\bu\in W^{1,0}_r(Q;\R^N)$ for which there exist vector functions  $\bg\in L^r(Q;\R^N)$ and $\bF\in L^r(Q;\M^{N\times n})$ such that
$$
\bu_t=\div \bF -\bg \quad \text{  a.e.  in } Q
$$
in the sense of distributions, that is, for each vector function $\chi\in C_0^\infty(Q)$ with $\chi(x,T)=0$
 holds
\begin{equation}\label{Fg}
\int_Q \bu\cdot\chi_t\, dxdt =\int_Q(\bF\cdot D\chi + \bg\cdot\chi)\, dxdt\,.
\end{equation}
The space $W^{1,r}_\ast(Q;\R^N)$  is endowed by the norm
$$
\|\bu\|_{W_\ast^{1,r}(Q)}=\|\bu\|_{W_r^{1,0}(Q)}+ \inf\left\{\left(\int_Q|\bF|^r+|\bg|^r  \right)^{1/r}\right\}
$$
where the infimum is taken over all $\bF$ and $\bg$ satisfying \eqref{Fg}.
The closure of $C_0^\infty(Q)$ with respect to  this norm is denoted by 
$\overset{\circ}{W}{}^{1,r}_\ast(Q;\R^N).$\\
3. $V_2(Q;\R^N)$ stands for the Banach space of all functions $\bu\in W^{1,0}_2(Q;\R^N)$ for which
$$
\|\bu\|_{V_2(Q;\R^N)}=\underset{t\in[0,T]}{\text{ess\,sup}}\,
 \|\bu(\cdot,t) \|_{2,\Omega} + \|D\bu\|_{2,Q}<\infty\,.
$$
4. $V^{1,0}_2(Q;\R^N)$ consists of all $\bu\in V_2(Q;\R^N)$ that are continuous in $t$ with respect to the norm of $L^2(\Omega;\R^N)$ 
$$
\lim_{\Delta t\to 0}\|\bu(\cdot, t+\Delta t) -\bu(\cdot,t)  \|_{2,\Omega}=0.
$$
The norm in $V^{1,0}_2(Q;\R^N)$ is given by 
$$
\|\bu\|_{V^{1,0}_2(Q;\R^N)}=\max_{t\in[0,T]}\|\bu(\cdot,t)\|_{2,\Omega}+\|D\bu\|_{2,Q}.
$$
 We consider the  Cauchy-Dirichlet problem for the strongly parabolic quasilinear system
\begin{equation}\label{CDP}
\begin{cases}
u_t^i-D_\alpha(A^{\alpha\beta}_{ij}D_\beta u^j+a^\alpha_i(x,t,\bu))=b_i(x,t,\bu,D\bu)&\text{a.a. }  (x,t)\in Q\\
\bu(x,t)=0&(x,t)\in \partial Q.
\end{cases}
\end{equation}
The principal coefficients satisfy  $A^{\alpha\beta}_{ij}\in L^\infty(Q)$ and
\begin{equation}\label{parabolic}
A^{\alpha\beta}_{ij}\xi^i_\alpha \xi^j_\beta\geq \nu |\xi|^2 \qquad \forall\  \xi\in \M^{N\times n},\quad \nu=\const>0. 
\end{equation}
In addition we suppose  that $A^{\alpha\beta}_{ij}$ are $(\delta,R)$-vanishing, that is 
\begin{equation}\label{deltaR}
\sup_{0<\rho\leq R}\sup_{(x,t)\in Q}\Xint-_{Q_\rho(x,t)} |A^{\alpha\beta}_{ij}(y,\tau)-\ol{A^{\alpha\beta}_{ij}}_{Q_\rho(x,t)}|^2dyd\tau\leq \delta^2
\end{equation}
which implies small $BMO$ norm of each $A^{\alpha\beta}_{ij}.$

The functions $a^\alpha_i(x,t,\bu),$ $b_i(x,t,\bu,\bz)$ are   Carath\'eodory ones   and verify the {\it controlled growth   conditions}   (see \cite{A,  LSU})
\begin{align}\label{growth1}
&\sum_{\underset{i\leq N}{\alpha \leq n}}|a_i^\alpha(x,t,\bu)|\leq \Lambda (\psi_1(x,t) +|\bu|^{\frac{n+2}{n}})\\ \label{growth2}
&\sum_{i\leq N} |b_i(x,t,\bu,\bz)|\leq \Lambda( \psi_2(x,t)+ |\bu|^{\frac{n+4}{n}}+|\bz|^{\frac{n+4}{n+2}}  )\,. 
\end{align}
with $ \psi_1\in L^p(Q), p>2 $ and  $\psi_2\in L^q(Q), q>\frac{2(n+2)}{n+4}.$

A vector function  $\bu \in \overset{\circ}{W}{}_\ast^{1,2}\cap L^\infty(Q;\R^N)$   is a weak solution to \eqref{CDP}  if for any    function $\chi \in \overset{\circ}{W}{}_\ast^{1,2}(Q;\R^N),$  $\chi(x,T)=0$  we have
\begin{align*}
\int_Qu^i(x,t)\chi_t(x,t)dxdt
&-\int_Q ( A^{\alpha\beta}_{ij}(x,t)D_\beta u^j(x,t)+a^\alpha_i(x,t,\bu))D_\alpha\chi^i(x,t)
dxdt\\
&+\int_Q b_i(x,t,\bu,D\bu)\chi^i(x,t) dxdt=0\,.
\end{align*}
\begin{thm}\label{th2}
Let $\bu\in \overset{\circ}{W}{}_\ast^{1,2}\cap L^\infty(Q;\R^N)$ be a weak solution to  \eqref{CDP} under the conditions \eqref{parabolic}-\eqref{growth2}. 
Then there exists a small positive constant  $\delta_0<1$ such that if $\Omega$ is $(\delta,R)$-Reifenberg flat
 and $A^{\alpha\beta}_{ij}$ are $(\delta,R)$-vanishing with
$0<\delta<\delta_0$ then  
\begin{align*}
&\bu\in   \overset{\circ}{W}{}_\ast^{1,r}\cap L^\infty(Q;\R^N) \qquad \text{ with }   r=\min\{p,q^{\ast\ast}\}\\
&\text{where }\quad q^{\ast\ast}=\begin{cases}
\frac{q(n+2)}{n+2-q} & \text{ if } q<n+2\\
\text{arbitrary large number }>1 & \text{ if } q\geq n+2.
\end{cases}
\end{align*}
 \end{thm}
\begin{pf}
The higher integrability of the gradient of the solution follows by the modification of the Gehring lemma due to Arkhipova \cite[Theorem~1]{A3}  which is very efficient  for the study of parabolic systems with controlled growth conditions in domains with boundary $\partial\Omega$ satisfying a kind of \eqref{A}-property. Recently,    similar result is obtained in \cite{DK} for domains having strongly Lipschitz boundary.  Since the Reifenberg flat domain can be approximated uniformly with Lipschitz domains with a small Lipschitz constant \cite[Lemma~5.1]{BW}) we have that  there exists $r_0>2$ such that for any solution 
$\bu\in   \overset{\circ}{W}{}_\ast^{1,2}(Q;\R^N)$  holds
$$
\|D\bu\|_{r,Q}\leq N \qquad \forall r\in[2,r_0)
$$
where $r_0$ and $N$ depend on the data of the problem and $\|D\bu\|_{2,Q}.$ Take a solution  $\bu\in \overset{\circ}{W}{}_\ast^{1,2}\cap L^\infty(Q;\R^N) $
of \eqref{CDP} and   fix  it in  the lower order terms.  Thus  we get the linearized   problem
\begin{equation}\label{linear2}
\begin{cases}
u^i_t-D_\alpha(A^{\alpha\beta}_{ij}D_\beta u^j )=f_i(x,t) +\div\, \A_i(x,t) & \text{ a.a. } (x,t)\in Q\\
\bu(x,t)=0 & \text{ on } \partial Q
\end{cases}
\end{equation}
where 
\begin{align*}
&\A_i(x,t)= (a_i^1(x,t,\bu),\ldots, a_i^n(x,t,\bu)),\quad
 \bA(x,t)=\{a_i^\alpha(x,t,\bu(x,t))\}_{i\leq N}^{\alpha\leq n}\\
&f_i(x,t)= b_i(x,t,\bu,D\bu),\quad \ff=(f_1(x,t),\ldots,f_N(x,t)).
\end{align*}
Making use of the conditions \eqref{growth1} and \eqref{growth2} we get
\begin{align*}
&\sum_{i\leq N} \|\A_i\|_{p,Q}\leq C\big(\|\psi_1\|_{p,Q}+\|\bu\|_{\infty,Q}^{\frac{n+2}{n}} \big)\\
&\sum_{i\leq N}\|f_i\|_{q_1,Q}\leq C\big(\|\psi_2\|_{q_1,Q} +\|\bu\|_{\infty,Q}^{\frac{n+4}{n}}+
\|D\bu\|_{\frac{q_1(n+4)}{n+2}}^{\frac{n+4}{n+2}}   \big)
\end{align*}
with   $p>2,$ $q_1=\min \big\{q,\frac{r_0(n+2)}{n+4} \big\}.$
Let $\Omega'\subset\R^n$ be $C^2$-domain such that $\Omega\Subset\Omega'$ and consider the cylinder 
$Q'=\Omega'\times(0,T).$ Suppose that $f_i(x,t)$ is extended as zero  out of $Q.$ 
It is well known (cf. \cite{LSU}) that for each $f_i\in L^{q_1}(Q')$ the linear problem
$$
\begin{cases}
F^i_t -\Delta F^i=F^i_t -D_\alpha(\delta^{\alpha\beta}D_\beta F^i)  = f_i(x,t) & \text{ a.a. } (x,t)\in Q'\\
F^i(x,t)=0 & \text{ on } \partial Q'
\end{cases}
$$
has a unique solution $F^i\in \overset{\circ}{W}{}^{2,1}_{q_1}(Q')$ and the estimate holds
$$
\|F^i\|_{q_1,Q}\leq \|F^i\|_{ \overset{\circ}{W}{}^{2,1}_{q_1}(Q')}\leq C\|f_i\|_{q_1,Q}.
$$
Denote by   $\bF=(F^1,\ldots,F^N),$  $\F_i^\alpha=(A^{\alpha\beta}_{ij} -\delta^{\alpha\beta}\delta_{ij})D_\beta F^j$ and 
 $\F_i=(\F_i^1,\ldots,\F_i^n).$ 
 Since the Sobolev trace theorem holds for domain  with Reifenberg flat boundary  we get that $F^i\vert_{S}\in
 W^{2-\frac1{q_1},1-\frac1{2q_1}}(S)$ where $S=\partial \Omega\times(0,T).$ Consider  the linear system
\begin{equation}\label{linear3}
\begin{cases}
(u^i-F^i)_t-D_\alpha(A^{\alpha\beta}_{ij}D_\beta(u^j-F^j) )&\\
\qquad =D_\alpha((A^{\alpha\beta}_{ij} -\delta^{\alpha\beta}\delta_{ij})D_\beta F^j )+\div\,(\A_i)\quad &\\
\qquad    =\div\,(\F_i+\A_i)  & \text{ a.a. } (x,t)\in Q\\
\bu-\bF =-\bF & \text{ on } \partial Q.
\end{cases}
\end{equation}
By the  imbedding theorems, $DF^i\in L^{q_1^{\ast\ast}}$  
 (cf. \cite[Ch.II,Lemma~3.3]{LSU}) and 
$$
\|DF^i\|_{q_1^{\ast\ast},Q}\leq C \|f^i\|_{q_1,Q}
$$
Hence $\F_i+\A_i\in L^{r_1}(Q)$ with $r_1=\min\{p,q_1^{\ast\ast}\}.$  Applying  \cite[Corollary~2.10]{BR} on the linear problem \eqref{linear3} we get
$$
\|D(\bu-\bF)\|_{r_1,Q}\leq C\big(1+\|D\bF\|_{q_1^{**},Q} + \|\bA\|_{p,Q} \big)\leq C\big(1+\|\ff\|_{q_1,Q}+\|\bA\|_{p,Q}\big).
$$
Hence 
$$
\|D\bu\|_{r_1,Q}\leq C(1+\|\ff\|_{q_1,Q}+\|\bA\|_{p,Q}).
$$
Applying  the bootstrapping arguments , we obtain as in the elliptic case that the theorem holds after $k$ iterations with 
$$
k\geq\min\left\{ \big[\frac{\log (r_0/(r_0-2))}{\log((n+4)/(n+2))} \big], \big[\log\frac{r_0[q(n+4)-2(n+2)]}{q(r_0-2)(n+4)} \big/ \log\frac{n+4}{n+2} \big]+1\right\}.
$$ 
\end{pf}

\subsection*{Acknowledgments.} The author is indebted to the referee for the valuable remarks which have lead to improvement of the article.

\end{document}